\documentclass[11pt,bezier,epsf]{amsart}
\usepackage{amsmath,amssymb,amsfonts}
\usepackage{graphicx}
\usepackage{cite}
\usepackage{subfigure}
\newbox\subfigbox
\makeatletter
{\def\caption##1{\gdef\subcapsave{\relax##1}}%
  \let\subcapsave\@empty%
  \setbox\subfigbox\hbox%
  \bgroup}%
{\egroup%
  \subfigure[\subcapsave]{\box\subfigbox}}%
\makeatother

 \newtheorem{thm}{Theorem}[section]

\theoremstyle{definition}
\newtheorem{example}[thm]{Example}
 \newtheorem{defn}[thm]{Definition}
\theoremstyle{remark}
 \newtheorem{rem}[thm]{Remark}
\numberwithin{equation}{section}
\begin{document}

\title[wrong intuitions]{how wrong intuitions about weak* topology and completions of a normed space pose serious problems}

\author[F. Naderi]{Fouad Naderi}

\address{Department of Mathematical and Statistical Sciences, University of
Alberta, Edmonton, Alberta, T6G 2G1, Canada.}
\email{naderi@ualberta.ca}


\begin{abstract}
We learn mathematics subjectively and must apply it objectively. But sometimes, we apply it subjectively by using wrong intuitions which may be elusive to our eyes. The aim of this note is to disclose the secretes of two kinds of these false intuitions and the opportunities they may provide. We first discuss the wrong assumption which says that each topology is uniquely determined by studying a very bad phenomenon happening in dealing with weak* topology. Then, we consider the problem of completing a normed space under two comparable norms, one being smaller than the other. Here, we show that the common belief contending that smaller norms give rise to larger completions is wrong. We then pose some serious questions arising from these wrong intuitions. As we will see finding fallacies are as important as major mathematical activities like proving and disproving.
\end{abstract}

\maketitle
%%%%%%%%%%%%%%%%%%%%%%%%%%%%%%%%%%%%%%%%%%%%%%%%%%%%%%%%%%%%%%%%%%%%%%%%%%%%%%%%%%%%%%%%%%%%%
\section{A serious problem with weak* topology}
We begin with reviewing some materials from \cite{Eymard} and \cite{Dixmier} which we need in the sequel.

Let $G$ be a locally compact group with a fixed left Haar measure $\mu$. Let
$L^{1}(G)$ be the group algebra of $G$ with convolution product. We define
$C^{\ast}(G)$, the group $C^{\ast}$-algebra of $G$, to be the completion of
$L^{1}(G)$ with respect to the norm
 $$\|f\|_{*}=sup\|\pi_{f}\|$$
where the supremum is taken over all nondegenerate $\ast$-representation $\pi$
of $L^{1}(G)$ as a $\ast$-algebra of bounded operators on a Hilbert space. One of such 
representations is called the left regular representation and is denoted by $\rho$. If
$\mathcal{B}(L^{2}(G))$ denotes the set of all bounded operators on the Hilbert space
$L^{2}(G)$, for every $f\in L^{1}(G)$, $\rho(f)$ is the bounded operator in $\mathcal{B}(L^{2}(G))$
defined by $\rho(f)(h)=f\ast h$, the convolution of $f$ and $h$ in $L^{2}(G)$.
The reduced group C*-algebra ${{C^{*}}}_{\rho}(G)$ is the completion of $L^{1}(G)$ 
with respect to $\|f\|_{\rho}=\|\rho(f)\|$. When $G$
is amenable (e.g., abelian and compact groups), one has ${C^{*}}_{\rho}
(G)=C^{*}(G)$. Denote the set of all continuous positive 
definite function on $G$ by $P(G)$, and
the set of all continuous function on $G$ with compact support by $C_{c}(G)$. We
define the Fourier-Stieltjes  algebra $B(G)$ of $G$ to be the
linear span of $P(G)$. Then, $B(G)$ is a Banach algebra with pointwise
multiplication and the norm of each $\phi \in B(G)$ is defined by
$$\|\phi\|=sup\{ | \int f(t)\phi(t) d\mu(t) | : f\in L^{1}(G), \|f\|_{*}\leq 1
\}.$$
In fact, $B(G)$ is the set of all matrix elements of all
continuous unitary representations of $G$. It is known that $B(G)=C^{\ast}(G)^{\ast}$, where the duality is given
by $\langle f,\phi\rangle = \int f(t)\phi(t) d\mu(t), f\in L^{1}(G), \phi \in
B(G)$. The closure of $\rho(L^{1}(G))$ 
in the weak operator topology of $\mathcal{B}(L^{2}(G))$ is called the von Neumann algebra generated by $\rho$
 and is denoted by $VN(G)$. For any other unitary representation $\pi$, the von Neumann algebra generated by $\pi$
is the closure of $\pi(L^{1}(G))$ 
in the weak operator topology of $\mathcal{B}({\mathcal{H}}_{\pi})$, where ${\mathcal{H}}_{\pi}$ is the Hilbert space 
of the representation $\pi$.

A von Neumann algebra $\mathcal{M}$ is called {\it atomic} if every nonzero
(self-adjoint) projection in $\mathcal{M}$ dominates a minimal nonzero projection
in $\mathcal{M}$. A locally compact group $G$ is called an
[AU]-group if the von Neumann algebra generated by any continuous unitary
representation of $G$ is atomic, i.e., every continuous unitary representation
of $G$ is completely decomposable. $G$ is called an [AR]-group if the von Neumann
algebra generated by $\rho$, $VN(G)$, is atomic. Evidently, every [AU]-group is
an [AR]-group (see \cite{Taylor} for more
details).

The set of all unitarily equivalent classes of
all irreducible representations of a C*-algebra $\mathcal{A}$ is called the
spectrum of $\mathcal{A}$ and is denoted by $\hat{\mathcal{A}}$. The dual
of a locally compact group $G$, denoted by $\hat{G}$, is the spectrum of 
the group C*-algebra $C^{*}(G)$. Similarly, the reduced
dual of a locally compact group $G$, denoted by $\hat{G}_{\rho}$, is the
spectrum of the reduced C*-algebra ${{C^{*}}}_{\rho}(G)$. The reduced dual has
many equivalent definitions (see \cite{Dixmier}).

A dual Banach space $E$ is said to have the {\it  uniform weak* Kadec-Klee property (UKK*)} if
for every $\epsilon >0$ there is $0<\delta<1$ such that for any subset $C$ of its closed unit ball
and any sequence $(x_n)$ in $C$ with $sep(x_n):=inf\{\rVert\ {x_n} -{ x_m}\rVert: n\neq m\}>\epsilon$, there is an $x$ in the weak* closure of $C$ with $\rVert x \rVert<\delta$. $E$ is said to have the {\it weak* Kadec-Klee property (KK*)} if the weak* and norm convergence for sequences coincide on the unit sphere of $E$. It is well known that UKK* implies KK* ( see \cite{Lennard} and \cite{Lau-Mah 88} for more details).

Now we are ready to give our first {\bf wrong result} based on a wrong intuition which we discuss after the theorem.

\begin{thm} \label{cpt}
{\bf (Wrong Theorem)} If $G$ is a second countable locally
compact [AU]-group, then $G$ is compact.
\end{thm}

{\bf Proof with one subtle gap.} At first, note that the C*-algebra $C^{*}(G)$ is separable since $G$ is second countable.  The dual Banach space $B(G)={C^{*}(G)}^{*}$ is the predual
${\mathcal{M}}_{*}$ of the von Neumann algebra $\mathcal{M}={B(G)}^{*}={C^{*}(G)}^{**}$. By the proof of \cite[Theorem 4.1]{Randrian}, one can embed ${\mathcal{M}}_{*}$ isometrically isomorphic into a sub-space of trace class operators $\mathcal{T}(\mathcal{H})$ for a suitable Hilbert space $\mathcal{H}$. Hence, by Lennard's theorem \cite[Theorem 2.4]{Lennard}, ${\mathcal{M}}_{*}$ is UKK* and KK*. On the other hand, ${\mathcal{M}}_{*}$ is the dual Banach space of the separable C*-algebra $C^{*}(G)$. Hence, the closed unit ball of 
${\mathcal{M}}_{*}$ is metrizable in the weak* topology \cite[p.134]{Conway}. So, the weak* topology of the closed unit ball is completely described by sequences. Since ${\mathcal{M}}_{*}$ is KK* , we deduce that the weak* and norm topology on the unit sphere of ${\mathcal{M}}_{*}=B(G)$ agree. Now, by \cite[Theorem 4.2]{Miao}, $G$ must be compact. $\blacksquare$
\begin{rem}
Theorem \ref{cpt} is of course wrong since the so called {\bf Fell group} is a counter example to its result (see \cite{Baggett} and \cite{Taylor}). But, we miss out one subtle remark if we do not examine its proof and content ourselves with just giving a counter example. There is only one gap in the end of the proof which makes the result false! In fact, the last part of the proof uses the erroneous assumption that weak*
topology is unique! At first, we endow ${\mathcal{M}}_{*}=B(G)$ with the relative weak* topology from $\mathcal{T}(\mathcal{H})$, then we shift to its weak* topology which comes from $C^{*}(G)$. But, the fact is that these weak* topologies may not be equal!, i.e., weak* topology is not uniquely determined (unlike weak topology). This is because preduals are {\bf NOT} unique (unlike duals)! A simple example which can explain the non-uniqueness of weak* topology  is $l^{1}$ with its two non-isomorphic preduals $c$ and $c_{0}$ ( $c$ is the space of all convergent sequences under the sup-norm topology and $c_{0}$ is its closed sub-space consisting of those elements converging to zero). On the other hand, von Neumann algebras (like $l^{\infty}$ ) have always unique preduals (and unique weak* topology) \cite[p.101]{Takesaki}.  Now, the question is why do we make such a mistake? At first, the notation $\sigma (A^{*}, A)$ for weak* topology is misleading since we know for any Banach space $A$ there is only one $A^{*}$. But, if we look at our symbol the other way round by putting $D=A^{*}$ and $P_{1}=A$, we now see that the weak* topology for $D$ is $\sigma (D, P_{1})$. The latter notation puts more emphasis on the dual space $D$ and asks which predual should be taken for it. In fact, we must interpret the notations $\sigma (A^{*}, A)$ or $\sigma (D, P_{1})$ as which predual ( and hence which weak* topology) is taken for the dual space. This is why the weak* topology is described by a pair of spaces; one is the space itself and the other is its predual. As we saw, every step in a proof must be taken on previously sturdy steps.
\end{rem}

Topics in \cite[p.101]{Takesaki} and the mistake in Theorem \ref{cpt} raise the following question:

{\bf Question 1.} Which Dual Banach spaces do have a unique predual? Is there any hope to characterize the uniqueness of predual in terms of algebraic and geometric properties of the dual space in question? 

It seems that the dual of a separable C*-algebras with discrete spectrum has a unique predual (!), and the proof of Theorem \ref{cpt} could be saved. For more details look at \cite[10.10.6]{Dixmier}, \cite[Theorem 3.4]{Baggett} and \cite[Theorem 3.1]{Naderi}.

%%%%%%%%%%%%%%%%%%%%%%%%%%%%%%%%%%%%%%%%%%%%%%%%%%%%%%%%%%%%%%%%%%%%%%%%%%%%%%%%%%%%%%%%%%%%% 
\section{Completing a normed space under two comparable norms}

Consider the following problem:

(*) Suppose that $X$ is a normed space under two norms ${\|.\|}_{s}$ and ${\|.\|}_{g}$ such that ${\|x\|}_{s}\leq {\|x\|}_{g}$ for any $x\in X$. Complete $X$ under these norms to obtain $X_s$ and $X_g$ (see \cite[p.82]{Rudin} for completion process). Now, a natural question is: which completion is larger? A widespread belief says  $X_s$ is larger since a smaller norm captures more Cauchy sequences. The following examples support this idea:

\begin{example} \label{completions}
(a) Endow $\mathbb{N}$ with its discrete topology and construct $c_{00}=C_{c}(\mathbb{N})$ as in \cite[p.132]{Folland}. $c_{00}$ is the set of all sequences with only finitely many non-zero terms. For any 
real number $1\leq p< \infty$ define ${\|.\|}_{p}$ as in \cite[Chapter 6]{Folland}. If one completes $c_{00}$ with respects to ${\|.\|}_{p}$, then (s)he obtains little $l^{p}$ \cite[p.217]{Folland}. For $1\leq p<q<\infty$, it is well-known that ${\|.\|}_{q}\leq {\|.\|}_{p}$ on $c_{00}$ \cite[p.186]{Folland}. The completion of $c_{00}$ with respect to ${\|.\|}_{q}$ and ${\|.\|}_{p}$ would be 
$l^{q}$ and $l^{p}$ respectively. In this case, ${\|.\|}_{q}\leq {\|.\|}_{p}$ and $l^{q}\supseteq l^{p}$.

(b) Let $X$ be a compact Hausdorff space with a Radon measure $\mu$ such that $\mu(X)=1$. By a similar discussion as in part (a) for $1\leq p<q<\infty$, we have ${\|.\|}_{p}\leq {\|.\|}_{q}$ and $L^{p}(\mu)\supseteq L^{q}(\mu)$ ( see \cite[Chapters 4, 6, and 7]{Folland} for more details).

(c) Let $e$ denotes the normalized Euclidean metric on $\mathbb{Q}$ and $d$ denotes the discrete metric on it. Obviously, we have $0\leq e(x,y)\leq d(x,y)\leq 1$ for any real $x$ and $y$. If we denote the completions of $\mathbb{Q}$ with respect to these metrics by ${\mathbb{Q}}_e$ and ${\mathbb{Q}}_d$, we see that ${\mathbb{Q}}_e=\mathbb{R}\supseteq \mathbb{Q}= {\mathbb{Q}}_d$.
\end{example}

But, the following example defies the intuition posed in (*):

\begin{example} \label{C*-compl}
The C*-algebras ${{C^{*}}}_{\rho}(G)$ and $C^{\ast}(G)$ in section 1 are completions of $L^{1}(G)$ under  $\|.\|_{\rho} \leq \|.\|_{*}$, but ${{C^{*}}}_{\rho}(G)$ is much {\bf smaller} than $C^{\ast}(G)$! In fact, the former is a quotient of the latter (see \cite{Eymard}). To know the reason, consider the inclusion map $i:(L^{1}(G),{\|.\|}_{*}) \hookrightarrow ({{C^{*}}}_{\rho}(G),{\|.\|}_{\rho})$ which, by th
the inequality $\|.\|_{\rho} \leq \|.\|_{*}$, can be extended to the continuous C*-algebraic homomorphism
$P:(C^{\ast}(G),{\|.\|}_{*}) \longrightarrow ({{C^{*}}}_{\rho}(G),{\|.\|}_{\rho})$. On the other hand, the continuous homomorphic image of C*-algebras are always closed \cite[p.81]{Murphy} and $L^{1}(G)$ is dense in ${{C^{*}}}_{\rho}(G)$, so $P$ is onto. Hence, ${{C^{*}}}_{\rho}(G)\cong \frac{C^{\ast}(G)}{Ker(P)}$.

\end{example}

\begin{rem}\label{clarify}
At this point, it is better to give the reason for wrong intuition in (*). We assumed  that
a smaller norm captures more Cauchy sequences, so the completion under smaller norm gives a bigger space. But, this is false since the number of Cauchy sequences is not the whole story in completion process! In fact, the number of different Cauchy
classes under the equivalence relation \cite[p.82]{Rudin} determines the size of completion. That is, a smaller norm may produce more zero Cauchy classes or makes more Cauchy sequences equivalent. Even, two non-equivalent Cauchy sequences in greater norm may be equivalent under smaller norm. In fact, in Examples \ref{completions} and \ref{C*-compl}, one uses indirect methods to decide which space is larger
(see the references provided) while the cardinality argument is futile. The best we can hope for is: under the conditions in (*), we always know there is a continuous homomorphism  $P:X_g \longrightarrow X_s$ extending the continuous inclusion map $i:(X,{\|.\|}_{g}) \hookrightarrow (X_{s},{\|.\|}_{s})$, but in general we do not know if such a $P$ is onto or not.
So, it is good to prove what seems intuitive and intuit what is proven.
\end{rem}
\begin{defn}
Under the assumptions in (*), we call the completions $X_s$ and $X_g$ order-preserving if they preserve the order of the inequality between norms; otherwise, we call them order-reversing.
\end{defn}
According to Remark \ref{clarify}, C*-completions of C*-norms always give rise to order-preserving completions. By Example
\ref{C*-compl}, there are order-reversing completions. Now, we ask:

{\bf Question 2.} Is there any rule to decide which completions are order-preserving and which are order-reversing?

{\bf Acknowledgments.} We would like to thank Professor Gero Fendler and Professor Eberhard Kaniuth for many 
fruitful correspondence.

%%%%%%%%%%%%%%%%%%%%%%%%%%%%%%%%%%%%%%%%%%%%%%%%%%%%%%%%%%%%%%%%%%%%%%%%%%%%%%%%%%%%%%%%%%%%%%%

%%%%%%%%%%%%%%%%%%%%%%%%%%%%%%%%%%%%%%%%%%%%%%%%%%%%%%%%%%%%%%%%%%%%%%%%%%%%%%%%%%%%%%%%%%%%%%%%%%%%%%%%%%%%%%%%%%%%%%%%%%%
\end{document}